\documentclass[12pt]{article}
\usepackage{latexsym,amsfonts,amssymb}
\usepackage{indentfirst}

\def\y{\begin{eqnarray*}}
\def\bd{\begin{description}}
\def\ey{\end{eqnarray*}}
\def\ebd{\end{description}}

\def\R{\mathbb{R}}
\def\N{\mathbb{N}}
\def\Z{\mathbb{Z}}

\begin{document}

\setlength{\baselineskip}{19pt}

\title{Existence and Asymptotic Stability of Periodic
Solutions for Impulsive Delay Evolution Equations
\thanks{Research supported by NNSF of China (11261053,11501455).
}}
\vskip3mm
\author{Qiang Li \footnote{Corresponding author. E-mail: lznwnuliqiang@126.com(Q.Li), nwnuweimei@126.com(M.Wei).}
,  Mei Wei
}

\date{\small
Department of Mathematics, Shanxi Normal University, \\
Linfen 041000, Peoples's Republic of China
 }

\maketitle

\begin{abstract}
\setlength{\baselineskip}{17pt}
In this paper, we are devoted to consider the periodic problem for the impulsive evolution equations with delay in Banach space. By using operator semigroups theory and fixed point theorem, we establish some new existence theorems of periodic mild solutions for the equations. In addition, with the aid of an integral inequality with impulsive and delay, we present essential conditions on the nonlinear and impulse functions to guarantee that the equations have an asymptotically stable $\omega$-periodic mild solution.

\vspace{8pt}

\vskip3mm
\noindent {\bf  MR(2010) Subject Classification:\ } 34G20; 34K30; 47H07; 47H08
\vskip5mm
\end{abstract}

\section{Introduction}

Let $X$ be a real Banach space with norm $\|\cdot\|$, $\mathcal{L}(X)$ stand for the Banach space of all bounded linear
 operators from $X$ to $X$ equipped
with its natural topology. Let $r>0$ be a constant, we  denote $PC([-r, 0], X)$ as the Banach space of piecewise continuous functions from $[-r, 0]$ to $X$
with finite points of discontinuity where functions are left continuous and have the right limits, with the sup-norm $\|\phi\|_{Pr}=\sup_{s\in [-r,0]}\|\phi(s)\|$.

In this paper, we consider the periodic problem  for the impulsive delay evolution equation  in  Banach space $X$
$$\left\{\begin{array}{ll}
u'(t)+Au(t)=F(t,u(t),u_{t}),\ t\in\R,\ t\neq t_{i},\\[8pt]
 \Delta u(t_{i})=I_{i}(u(t_{i})),\ \qquad i\in\Z,
 \end{array} \right.\eqno(1.1)$$
where  $A:D(A)\subset X\to X$ is a closed linear operator and $-A$ generates a $C_0$-semigroup $T(t)(t\geq0)$  in $X$; $F:\R\times X\times PC([-r, 0], X)\rightarrow X$ is a nonlinear mapping which is $\omega$-periodic in $t$; $u_{t}\in PC([-r,0],X)$ is the history function defined by  $u_{t}(s)=u(t+s)$ for $s\in [-r,0]$; $p\in \N$ denotes the number of impulsive points between $0$ and $\omega$, $0<t_{1}<t_{2}<\cdots<t_{p}<\omega<t_{p+1}$ are given numbers satisfying $t_{p+i}=t_{i}+\omega (i\in\Z)$; $\Delta u(t_{i})=u(t_{i}^{+})-u(t_{i}^{-})$ represents the jump of the function $u$ at $t_{i}$, $u(t_{i}^{+})$ and $u(t_{i}^{-})$ are the right and left limits of $u(t)$ at $t_{i}$, respectively; $I_i : X\rightarrow X(i\in \Z)$ are continuous functions satisfying $I_{i+p}=I_{i}$.

The theory of partial differential equations with delays has extensive physical background and realistic mathematical model, and it has undergone a rapid development in the last fifty years.
 The evolution equations with delay are more realistic than the equations without delay
 in describing numerous phenomena observed in nature,
  hence the numerous properties of their solutions have been studied,
  see \cite{Hale1993,Wu1996} and references therein for more comments.

One of the important research directions related to the asymptotic behavior of the solutions for the evolution equations with delay is to find conditions for the existence and stability of  $\omega$-periodic solutions in the case that the nonlinear mapping is $\omega$-periodic function in $t$. In the last few decades, the existence and asymptotic stability of periodic solutions have been investigated by some authors(see \cite{Burton1991,Xiang1992,Liu1998,Liu2000,Liu2003, Li2011,LiLi2016,Huy} and the references therein).  In \cite{Burton1991}, under the assumption that the solutions of the associated homotopy equations were uniformly bounded, Burton and Zhang obtained the existence of periodic solutions of an abstract delay evolution equation. In \cite{Xiang1992}, Xiang and Ahmed showed an existence result of periodic solution to the delay evolution equations in Banach spaces under the assumption that the corresponding initial value problem had a priori estimate. In \cite{Liu1998,Liu2000,Liu2003}, Liu studied periodic solutions by using bounded solutions or ultimate bounded solutions for  delay evolution equations in Banach spaces. In \cite{Huy}, Huy and Dang  studied the existence, uniqueness and stability of periodic solutions to a partial functional differential equation in Banach space in the case that the nonlinear function satisfied Lipschitz-type condition. Specially, in \cite{Li2011}, Li discussed the existence and asymptotic stability of periodic solutions to the evolution equation with multiple delays in a Hilbert space. By using the analytic semigroups theory and the integral inequality with delays, the author obtained the essential conditions
on the nonlinearity $F$ to guarantee that the equation has $\omega$-periodic solutions or an asymptotically stable $\omega$-periodic solution.

On the other hand, it is well known that impulsive evolution equation has
an extensive physical, chemical, biological, engineering background and realistic mathematical model, and hence has been emerging as an important area of investigation in the last few decades. Since the end of last century, theory of impulsive evolution equation in Banach space
has been largely developed(see \cite{Guo1993,Ahmed01,Ahmed2006,Yang2006,Liang2009,Liang2011} and  the references therein).
We would like to mention that Liang, et al. \cite{Liang2011} studied the periodic solutions to a kind of impulsive evolution equation with delay in Banach spaces.
In the case that the nonlinear function satisfied
 Lipschitz conditions, the authors obtained that the evolution equation has periodic solution by using the ultimate boundedness of solutions and Horn's fixed point theorem. Recently, in \cite{Liang2017A,Liang2017D},
 Liang et al.  studied  nonautonomous  evolutionary
 equations with time delay and impulsive. Under the nonlinear term
 satisfying continuous and Lipschitzian, they proved the existence theorem for periodic mild solutions to the nonautonomous delay evolution equations.

Although there have been many meaningful results on the delay or impulsive
 evolution equation periodic problem in Banach space, to our knowledge, these results have  relatively large limitation. First of all, the most popular approach is the use of ultimate boundedness of solutions and the compactness of Poincar\'{e} map realized through some compact
embeddings.  However, in some concrete applications, it is difficult to choose an appropriate initial conditions to guarantee the boundedness of the solution.
Secondly, we observe that the most popular condition imposed on the nonlinear term $F$ is its Lipschitz-type condition. In fact, for equations arising in complicated reaction-diffusion processes, the nonlinear function $F$ represents the source of material or population, which dependents on time in diversified manners in many contexts. Thus, we may not hope to have the Lipschitz-type condition of $F$. Finally, there are few papers to study
the asymptotically stable of periodic solutions for the impulsive evolution equations with delay.

Motivated by the papers mentioned above, we consider the periodic problem for impulsive delay evolution equation (1.1) in Banach space. By using periodic extension and fixed point theorem, we study the existence of $\omega$-periodic mild solutions for Eq.(1.1). It is worth mentioning that the assumption of prior boundedness of solutions is not employed and the nonlinear term $F$  satisfies some growth condition, which is weaker than Lipschitz-type condition. On the other hand, by means of an integral inequality with
impulsive and delay, we present the asymptotic stability result for Eq.(1.1), which will make up the research in this area blank.

The rest of this paper is organized as follows. In Section 2, we collect some known definitions and notions, and then provide preliminary results which will be used throughout this paper. In Section 3, we apply the operator semigroup theory to find the $\omega$-periodic mild solutions for Eq.(1.1) and in Section 4, by strengthening the condition, we obtain the global asymptotic stability theorems for Eq.(1.1). In the last section, we give an example to illustrate the applicability of abstract results obtained in Section 3 and Section 4.

\section{Preliminaries}

Throughout this paper, we assume that $X$ is a Banach space with norm $\|\cdot\|$.

Now, we recall some notions and properties of operator semigroups, which are  essential for us. For the detailed theory of operator semigroups, we refer to \cite{Pazy}. Assume that $A:D(A)\subset X\rightarrow X$ is a closed linear operator and  $-A$ is the infinitesimal generator of a $C_{0}$-semigroup $T(t)(t\geq0)$ in  $X$. Then there exist $M>0$ and $\nu\in \R$ such that
$$\|T(t)\|\leq Me^{\nu t},\quad t\geq0.\eqno(2.1)$$
Let
$$\nu_{0}=\inf\{\nu\in \R |\ \mathrm{There\ exists}\ M>0\ \mathrm{ such \ that}\  \|T(t)\|\leq Me^{\nu t},\ \forall t\geq0\},$$
 then $\nu_{0}$ is called the growth exponent of the semigroup $T(t)(t\geq0)$. If $\nu_{0}<0$, then $T(t)(t\geq0)$ is called an exponentially
stable $C_0$-semigroup.

If $C_{0}$-semigroup $T(t)$ is continuous in the uniform operator topology for every $t>0$ in $X$, it is well known that $\nu_{0}$ can also be determined by $\sigma(A)$ (the resolvent set of $A$),
$$\nu_{0}=-\inf\{\mathrm{Re} \lambda |\ \lambda\in\sigma(A)\},\eqno(2.2)$$
where $-A$ is the infinitesimal generator of $C_{0}$-semigroup $T(t)(t\geq0)$. We know that $T(t)(t\geq0)$ is continuous in the uniform operator topology for $t>0$ if $T(t)(t\geq0)$ is compact semigroup (see \cite{Triggiani}).

Let $J$ denote the infinite interval $[0,+\infty)$, from \cite{Pazy}, it follows that when $x_{0}\in D(A)$ and $h\in C^{1}(J,X)$, the following initial value problem of the linear evolution equation
$$\left\{\begin{array}{ll}
u'(t)+Au(t)=h(t),\ t\in J,\\[8pt]
u(0)=x_{0}
 \end{array} \right.\eqno(2.3)$$
  has a unique classical solution
 $u\in C^{1}(J,X)\cap C(J,X_{1})$, which can be expressed by
 $$u(t)=T(t)x_{0}+\int^{t}_{0}T(t-s)h(s)ds, \eqno(2.4)$$
where $X_{1}=D(A)$ is
Banach space  with the graph norm $\|\cdot\|_{1} = \|\cdot\|+\|A\cdot\|$.
Generally, for $x_{0}\in X$ and $h\in C(J,X)$,
the function $u$ given by (2.4) belongs to $C(J,X)$ and it is
called a mild solution of the linear evolution equation (2.3).

Let $\widetilde{D}=\{t_{1},t_{2},\cdots,t_{p}\}\subset [0,\omega]$, where $p\in N$ denotes the number of impulsive points between $[0,\omega]$. Write
\begin{eqnarray*}PC([0,\omega],X):=\Big\{u:[0,\omega]\to X| u\ \mathrm{is\ continuous\ at}\ t\in [0,\omega]\setminus \widetilde{D}, \\[8pt]
u\ \mathrm{is\ continuous\ from\ left\ and\ has\ right\ limits\ at}\ t\in \widetilde{D}\Big\}\end{eqnarray*}
and
\begin{eqnarray*}PC_{\omega}(\R,X)&:=&\Big\{u:\R\to X| u\ \mathrm{is\ continuous\ at}\ t\in \R\setminus \{t_{i}\},u\ \mathrm{is\ continuous\ from}\\[8pt]
&& \mathrm{left\ and\ has\ right\ limits\ at}\ t_{i}, i\in \Z, u(t+\omega)=u(t)\ \mathrm{for}\ t\in \R \Big\}.\end{eqnarray*}
It is clear that the restriction of $PC_{\omega}(\R,X)$ on $[0,\omega]$ is $PC([0,\omega],X)$. Set $$\|u\|_{PC}=\max\Big\{\sup_{t\in [0,\omega]}\|u(t+0)\|,\ \sup_{t\in [0,\omega]}\|u(t+0)\|\Big\}.$$
It can be seen that equipped with the norm $\|\cdot\|_{PC}$, $PC_{\omega}(\R,X)$ (or $PC([0,\omega],X)$) is a Banach space.

Given $h\in PC_{\omega}(\R,X)$, we consider the existence of $\omega$-periodic mild solution for the linear impulsive evolution equation in $X$
$$\left\{\begin{array}{ll}
u'(t)+Au(t)=h(t),\qquad t\in\R,\ t\neq t_{i},\\[8pt]
\Delta u(t_{i})=v_{i},\  \qquad\qquad \quad i\in\Z,
 \end{array} \right.\eqno(2.5)$$
 where $v_{i}\in X$ satisfy
 $v_{i+p}=v_{i}(i\in\Z)$.

 \vskip3mm
\noindent\textbf{Lemma 2.1.} Let $-A$ generate an exponentially stable $C_{0}$-semigroup $T(t)(t\geq0)$ in $X$ and $\nu_{0}$ be a growth index of the semigroup $T(t)(t\geq0)$. Then the linear impuisive evolution equation (2.5) exists a unique $\omega$-periodic mild
solution $u:=Ph\in PC_{\omega}(\R,X)$.
Furthermore, the operator $P:PC_{\omega}(\R,X)\to PC_{\omega}(\R,X)$ is a bounded linear operator.

 \vskip2mm
\noindent\textbf{Proof.} Firstly, for $x_{0}\in X$, we consider the existence of mild solution for the initial value problem of the linear impulsive evolution equation
$$\left\{\begin{array}{ll}
u'(t)+Au(t)=h(t),\ \quad t\in\ J\setminus \{t_{1},t_{2},\cdots\},\\[8pt]
\Delta u(t_{i})=v_{i},\ \qquad i=1,2,\cdots,\\[8pt]
u(0)=x_{0}.
 \end{array} \right.\eqno(2.6)$$
Let $t_{0}=0, v_{0}=\theta$ and $J_{i}=(t_{i-1},t_{i}],\ i=1,2,\cdots$. If $u\in PC(J,X)$ is a mild solution of Eq.(2.6), then the restriction of $u$ on $J_{i}$ satisfies the initial value problem of linear evolution equation without impulse
$$\left\{\begin{array}{ll}
u'(t)+Au(t)=h(t),\ \quad t_{i-1}<t\leq t_{i},\\[8pt]
u(t_{i-1}^{+})=u(t_{i-1})+v_{i-1},\quad i=1,2,\cdots.
\end{array} \right.\eqno(2.7)$$
Thus, the  initial value problem (2.7) exists a unique mild solution $u$ on $J_{i}$ which can be expressed by
$$u(t)=T(t-t_{i-1})(u(t_{i-1})+v_{i-1})+\int^{t}_{t_{i-1}}T(t-s)h(s)ds,\ \ t\in J_{i} .\eqno(2.9)$$
Iterating successively in the above equation with $u(t_{j}),j=i-1,i-2,\cdots,1$, we can verify  that $u$ satisfies
$$u(t)=T(t)x_{0}+\int^{t}_{0}T(t-s)h(s)ds+\sum_{0\leq t_{i}<t}T(t-t^{+}_{i})v_{i},\qquad t\in J.\eqno(2.10)$$
Inversely, we can see that the function $u \in PC(J, X)$ defined by (2.9) is a mild solution of the initial value problem (2.6).

Secondly, we demonstrate that the linear impuisive evolution equation (2.5) exists a unique $\omega$-periodic mild solution. It is clear that the $\omega$-periodic mild solution of Eq.(2.5) restricted on $J$ is the mild solution of the initial value problem (2.6) with the initial value
$$x_{0}:=u(0)=u(\omega),$$
namely
$$(I-T(\omega))x_{0}=\int_{0}^{\omega}T(\omega-s)h(s)ds+\sum_{i=1}^{p}T(\omega-t^{+}_{i})v_{i}.\eqno(2.10)$$
For any $\nu\in(0,|\nu_{0}|)$, there exists $M>0$ such that
$$\|T(t)\|\leq Me^{-\nu t}\leq M, \quad t\geq0.\eqno(2.11)$$
In $X$, define an equivalent norm $|\cdot|$ by $|x|=\sup\limits_{0\leq s\leq t}\|e^{\nu t}T(t)x\|$,
then $\|x\|\leq |x|\leq M\|x\|$. By $|T(t)|$ we denote the norm of $T(t)$
in $(X,|\cdot|$), it is easy to obtain that $|T(t)|\leq e^{-\nu t}$ for $t\geq 0 $, which implies that $(I-T(\omega))$ has bounded inverse operator
$(I-T(\omega))^{-1}$, and
$$|(I-T(\omega))^{-1}|\leq \frac{1}{1-e^{-\nu\omega}}.\eqno(2.12)$$
Therefore, there exists a unique initial value $$x_{0}=(I-T(\omega))^{-1}\Big(\int^{\omega}_{0}T(\omega-s)h(s)ds
+\sum^{p}_{i=1}T(\omega-t^{+}_{i})v_{i}\Big):=B(h),\eqno(2.13)$$
 such that the mild solution $u$ of Eq.(2.6) given by (2.9) satisfies the periodic boundary condition $u(0)=u(\omega)=x_{0}$.

For $t\in J$, by the semigroup properties of $T(t)$, we have
\begin{eqnarray*}
u(t+\omega)&=&T(t+\omega)u(0)+\int_{0}^{t+\omega}T(t+\omega-s)h(s)ds+\sum_{0<t_{i}<t+\omega}T(t+\omega-t^{+}_{i})v_{i}\\
&=&T(t+\omega)u(0)+\int_{0}^{\omega}T(t+\omega-s)h(s)ds+\int_{\omega}^{t+\omega}T(t+\omega-s)h(s)ds\\
&\ &+\sum_{i=1}^{p}T(t+\omega-t^{+}_{i})v_{i}+\sum_{\omega<t_{i}<t+\omega}T(t+\omega-t^{+}_{i})v_{i}\\
&=&T(t)\Big(T(\omega)u(0)+\int_{0}^{\omega}T(\omega-s)h(s)ds+\sum_{i=1}^{p}T(\omega-t^{+}_{i})v_{i}\Big)\\
&\ &+\int_{0}^{t}T(t-s)h(s)ds+\sum_{0<t_{i}<t}T(t-t^{+}_{i})v_{i}\\
&=&T(t)u(0)+\int_{0}^{t}T(t-s)h(s)ds+\sum_{0<t_{i}<t}T(t-t^{+}_{i})v_{i}\\
&=&u(t).
\end{eqnarray*}
Therefore, the $\omega$-periodic extension of $u$ on $\R$ is a unique $\omega$-periodic mild solution of Eq.(2.5).

Finally, by (2.9) and (2.13), we obtain
\begin{eqnarray*}
u(t)&=&T(t)x_{0}+\int_{0}^{t}T(t-s)h(s)ds+\sum_{0<t_{i}<t}T(t-t^{+}_{i})v_{i}\\
&=&T(t)(I-T(\omega))^{-1}\Big(\int^{\omega}_{0}T(\omega-s)h(s)ds
+\sum^{p}_{i=1}T(\omega-t^{+}_{i})v_{i}\Big)\\
&&+\int_{0}^{t}T(t-s)h(s)ds+\sum_{0<t_{i}<t}T(t-t^{+}_{i})v_{i}\\
&=&(I-T(\omega))^{-1}\int_{t-\omega}^{t}T(t-s)h(s)ds\\
&&+(I-T(\omega))^{-1}\sum^{p}_{i=1}T(t+\omega-t_{i})v_{i}+\sum_{0<t_{i}<t}T(t-t_{i})v_{i}\\
&=&(I-T(\omega))^{-1}\int_{t-\omega}^{t}T(t-s)h(s)ds+
(I-T(\omega))^{-1}\sum_{t-\omega\leq t_{i}<t}T(t-t_{i})v_{i}\\
&:=&Ph(t),
\end{eqnarray*}
thus, it is  easy to prove that the solution operator $P: PC_{\omega}(\R,X)\to PC_{\omega}(\R,X)$ is a bounded linear operator.  This completes the proof of Theorem 2.1.  \hfill$\Box$

\section{ The Existence and Uniqueness}
In this section, we discuss the existence of $\omega$-periodic mild solution to Eq. (1.1) for the case that the semigroup $T(t)(t\geq0)$
generated by $-A$ is a  compact semigroup, which implies that $T(t)$ is a compact operator for any $t>0$.

Now, we are in a position to state and prove our main results of this section.
\vskip3mm
\noindent\textbf{Theorem 3.1.} \emph{Let $X$ be a Banach space, $-A$ generates an exponentially stable compact semigroup $T(t)(t\geq0)$ in $X$. Assume that  $F:\R\times X\times PC([-r, 0], X)\rightarrow X$  is continuous and $F(t,\cdot,\cdot)$ is $\omega$-periodic in $t$, $I_{k}\in C(X,X)(k\in \mathbb{Z})$ satisfies $I_{k+p}=I_{k}$, $p$ is the number of impulsive points between $[0,\omega]$. If the following condition is satisfied:
\vskip3mm
\noindent(H1) there exist nonnegative constants $c_{0},c_{1},c_{2}$
such that
$$\|F(t,x,\phi)\|\leq c_{0}+c_{1}\|x\|+c_{2}\|\phi\|_{Pr}, \ \ t\in\R,x\in X,\phi\in PC([-r,0],X),$$
\vskip3mm
\noindent(H2) for every $I_{k}$, $I_{k}(\theta)=\theta$,  and there exist positive constants $a_{k}$ such that
$$\|I_{k}(x)-I_{k}(y)\|\leq a_{k}\|x-y\|,\ \  a_{k+p}=a_{k}\ \ \ x,y\in X,k\in \Z,$$
\vskip3mm
\noindent(H3) $(c_{1}+c_{2})+\frac{1}{\omega}\sum\limits_{k=1}^{p}a_{k}<\frac{|\nu_{0}|}{M}$,
\vskip3mm
\noindent then Eq. (1.1) has at least one $\omega$-periodic
mild solution $u$.}
\vskip3mm
\noindent\textbf{Proof } Define an operator $Q:PC_{\omega}(\R,X)\to PC_{\omega}(\R,X)$ by
\begin{eqnarray*}
\qquad (Qu)(t)&=&(I-T(\omega))^{-1}\int_{t-\omega}^{t}T(t-s)f(s,u(s),u_{s})ds\\
&&+(I-T(\omega))^{-1}\sum_{t-\omega\leq t_{k}<t}T(t-t_{k})I_{k}(u(t_{k})),\ \ \qquad t\in\R. \qquad\ (3.1)\end{eqnarray*}
From Lemma 2.1 and the definition of $Q$, we easily obtain that the $\omega$-periodic mild solution of Eq.(1.1) is equivalent to the fixed point
of $Q$. In the following, we will prove $Q$ has a fixed point by applying the
 fixed point theorem.

It is not difficult to prove that $Q$ is continuous on $PC_{\omega}(\R,X)$. In fact, let $\{u_{n}\}\subset PC_{\omega}(\R,X) $ be a sequence
such that $u_{n} \to u\in PC_{\omega}(\R,X)$ as $n\to \infty$,  hence, for every $t\in\R$, we have $u_{n}(t) \to u(t)\in X$ and $u_{n,t}\to u_{t}\in PC([-r,0],X)$ as $n\to \infty$.
From
$F: \R\times X\times PC([-r,0],X)\to X$ is
continuous, and $I_{k}\in C(X,X)(k\in \mathbb{Z})$, it follows that
$$F(t,u_{n}(t),u_{n,t})\to f(t,u(t),u_{t})),\ \ \ n\to \infty,\eqno(3.2)$$
and
$$I_{k}(u_{n}(t_{k}))\to I_{k}(u(t_{k})),\ \ \ n\to \infty. \eqno(3.3)$$
By (3.1)-(3.3) and the Lebesgue dominated convergence theorem, for every $t \in\R$, we have
\begin{eqnarray*}
&&\|(Qu_{n})(t)-(Qu)(t)\|\\[8pt]
&\leq& \|(I-T(\omega))^{-1}\|\cdot\Big\|\sum_{t-\omega\leq t_{k}<t}T(t-t_{k})(I_{k}(u_{n}(t_{k}))-I_{k}(u(t_{k})))\Big\|\\
&&+\|(I-T(\omega))^{-1}\|\cdot\Big\|\int_{t-\omega}^{\omega}T(t-s)\cdot\Big(F(s,u_{n}(s),u_{n,s})ds
-F(s,u(s),u_{s})\Big)ds\Big\|\\
&\leq& C\Big(\sum_{t-\omega\leq t_{k}<t}\|T(t-t_{k})\|\cdot\|I_{k}(u_{n}(t_{k}))-I_{k}(u(t_{k}))\|\\
&&+\int_{t-\omega}^{t}\|T(t-s)\|\cdot\|F(s,u_{n}(s),u_{n,s})-F(s,u(s),u_{s})\|ds\Big)\\[8pt]
&\to& 0,\ \  (n\in\infty),
\end{eqnarray*}
 which implies that $Q:PC_{\omega}(\R,X)\to PC_{\omega}(\R,X)$ is  continuous,  where $C=\|(I-T(\omega))^{-1}\|$, by the proof of Lemma 2.1, one can obtain
$$\|(I-T(\omega))^{-1}\|\leq |(I-T(\omega))^{-1}|\leq\frac{1}{1-e^{\nu_{0}\omega}}.$$
 For any  $R>0$, let
$$\overline{\Omega}_{R}=\{u\in PC_{\omega}(\R,X)\ |\  \|u\|_{PC}\leq R\}.\eqno(3.4)$$
Note that $\overline{\Omega}_{R}$ is a closed ball in $PC_{\omega}(\R,X)$ with centre $\theta$ and radius $r$. Now, we  show that there is a positive constant $R$ such that $Q(\overline{\Omega}_{R})\subset\overline{\Omega}_{R}$. If this were not case, then for any $R>0$, there exist $u\in\overline{\Omega}_{R}$ and $t\in \R$ such that $\|(Qu)(t)\|>R$. Thus, we see by (H1) and (H2) that
\begin{eqnarray*}
R&<&\|(Qu)(t)\|\\[8pt]
&\leq&\|(I-T(\omega))^{-1}\|\cdot\int_{t-\omega}^{\omega}\|T(t-s)\|\cdot\|F(s,u(s),u_{s})\|ds\\
&&+\|(I-T(\omega))^{-1}\|\cdot\sum_{t-\omega\leq t_{k}<t}\|T(t-t_{k})\|\cdot\|I_{k}(u(t_{k}))\|\\[8pt]
&\leq&\frac{1}{1-e^{\nu_{0}\omega}}\cdot\int_{t-\omega}^{t}Me^{\nu_{0}(t-s)}(c_{0}+c_{1}\|u(s)\|+c_{2}\|u_{s}\|_{Pr})ds\\
&&+\frac{1}{1-e^{\nu_{0}\omega}}\cdot\sum_{k=1}^{p}M e^{\nu_{0}(t-t_{k})}a_{k}\|u(t_{k})\|\\[8pt]
&\leq&-\frac{M}{\nu_{0}}(c_{0}+(c_{1}+c_{2})\|u\|_{PC})-\frac{M}{\nu_{0}\omega}\sum_{k=1}^{p}a_{k}\|u\|_{PC}\\[8pt]
&\leq&-\frac{M}{\nu_{0}}(c_{0}+(c_{1}+c_{2})R)-\frac{M}{\nu_{0}\omega}\sum_{k=1}^{p}a_{k}R.
\end{eqnarray*}
Dividing on both sides by $R$ and taking the lower limit as $R\rightarrow\infty$, we have
$$(c_{1}+c_{2})+\frac{1}{\omega}\sum_{k=1}^{p}a_{k}\geq -\frac{\nu_{0}}{M},\eqno(3.5)$$
which contradicts (H3). Hence, there is a positive constant $R$ such that $Q(\overline{\Omega}_{R})\subset\overline{\Omega}_{R}$.

In order to show that the operator $Q$ has a fixed point on $\overline{\Omega}_{R}$,  we also
introduce the decomposition $Q=Q_{1}+Q_{2}$, where
$$Q_{1}u(t):=(I-T(\omega))^{-1}\int_{t-\omega}^{t}T(t-s)F(s,u(s),u_{s})ds,\ \eqno(3.6)$$
$$Q_{2}u(t):=(I-T(\omega))^{-1}\sum_{t-\omega\leq t_{k}<t}T(t-t_{k})I_{k}(u(t_{k})).\eqno(3.7)$$
Then we will prove that $Q_{1}$ is a compact operator and $Q_{2}$ is a contraction.

Firstly, we prove that $Q_{1}$ is a compact operator. Clearly,  $Q_{1}$ is continuous and
 $Q_{1}$ maps
 $\overline{\Omega}_{R}$ into a bounded set in $PC_{\omega}(\R,X)$. Now, we demonstrate that $ Q_{1}(\overline{\Omega}_{R})$ is equicontinuous. For every
$u\in\overline{\Omega}_{R}$, by the periodicity of $u$, we only consider it on $[0,\omega]$.
Set $ 0\leq t_{1}<t_{2}\leq\omega$, we get that
\begin{eqnarray*}
&&Q_{1}u(t_{2})-Q_{1}u(t_{1})\\[8pt]
&=&(I-T(\omega))^{-1}\int^{t_{2}}_{t_{2}-\omega}T(t_{2}-s)F(s,u(s),u_{s})ds\\
&~&-(I-T(\omega))^{-1}\int^{t_{1}}_{t_{1}-\omega}T(t_{1}-s)F(s,u(s),u_{s})ds\\[8pt]
&=&(I-T(\omega))^{-1}\int^{t_{1}}_{t_{2}-\omega}(T(t_{2}-s)-T(t_{1}-s))F(s,u(s),u_{s})ds\\
&~&-(I-T(\omega))^{-1}\int^{t_{2}-\omega}_{t_{1}-\omega}T(t_{1}-s)F(s,u(s),u_{s})ds\\[8pt]
&~&+(I-T(\omega))^{-1}\int^{t_{2}}_{t_{1}}T(t_{2}-s)F(s,u(s),u_{s})ds\\[8pt]
&:=&I_{1}+I_{2}+I_{3}.
\end{eqnarray*}
It is clear that
$$\|Q_{1}u(t_{2})-Q_{1}u(t_{1})\|\leq \|I_{1}\|+\|I_{2}\|+\|I_{3}\|.\eqno(3.8)$$
Thus, we only need to check $\|I_{i}\|$ tend to $0$ independently of $u\in\overline{\Omega}_{R}$
when $t_{2}-t_{1}\rightarrow 0,i=1,2,3$.

From the condition (H1), it follows that there is a constant $M'>0$ such that
$$\|F(t,u(t),u_{s})\|\leq M',\ \ u\in \overline{\Omega}_{R},\ t\in\R.$$
Combined this fact with the equicontinuity of the semigroup $T(t)(t\geq0)$, we have
\begin{eqnarray*}
\|I_{1}\|&\leq&C\cdot\int^{t_{1}}_{t_{2}-\omega}\|(T(t_{2}-s)-T(t_{1}-s))\|\cdot\|F(s,u(s),u_{s})\|ds\\[8pt]
&\leq & CM'\int^{t_{1}}_{t_{2}-\omega}\|(T(t_{2}-s)-T(t_{1}-s))\| ds\\[8pt]
&\rightarrow&0, \ \mathrm{as} \ t_{2}-t_{1}\rightarrow 0,\\[8pt]
\|I_{2}\|&\leq&C\cdot\int^{t_{2}-\omega}_{t_{1}-\omega}\|T(t_{1}-s)\|\cdot\|F(s,u(s),u_{s})\|ds\\[8pt]
&\leq &CMM'(t_{2}-t_{1})\\[8pt]
&\rightarrow&0,  \ \mathrm{as} \ t_{2}-t_{1}\rightarrow 0,\\[8pt]
\|I_{3}\|&\leq&C\cdot\int^{t_{2}}_{t_{1}}\|T(t_{2}-s)\|\cdot\|F(s,u(s),u_{s})\|ds\\[8pt]
&\leq &CMM'(t_{2}-t_{1})ds \\[8pt]
&\rightarrow&0,  \ \mathrm{as} \ t_{2}-t_{1}\rightarrow 0.
\end{eqnarray*}
Hence, $\| Q_{1}u(t_{2}) -Q_{1}u(t_{1})\|$ tends to $0$
 independently of $u\in \overline{\Omega}_{R}$
as $t_{2}- t_{1}\rightarrow0$,
which means that $Q_{1}(\overline{\Omega}_{R})$ is equicontinuous.

It remains to show that $(Q_{1}\overline{\Omega}_{R})(t)$ is relatively compact in $X$ for all $t\in \R$.
To do this, we define a set $(Q_{\varepsilon}\overline{\Omega}_{r})(t)$ by
$$(Q_{\varepsilon}\overline{\Omega}_{r})(t):=\{(Q_{\varepsilon}u)(t)|\ u\in \overline{\Omega}_{r},\ 0<\varepsilon<\omega,\ t\in \R\}, \eqno(3.9)$$
where
\begin{eqnarray*}
(Q_{\varepsilon}u)(t)&=&(I-T(\omega))^{-1}\int_{t-\omega}^{t-\varepsilon}T(t-s)F(s,u(s),u_{s})ds\\
&=&T(\varepsilon)(I-T(\omega))^{-1}\int_{t-\omega}^{t-\varepsilon}T(t-s-\varepsilon)F(s,u(s),u_{s})ds.
\end{eqnarray*}
 Since the operator $T(\varepsilon)$ is compact in $X$, thus, the set $(Q_{\varepsilon}\overline{\Omega}_{R})(t)$ is relatively compact in $X$.  For any $u\in \overline{\Omega}_{R}$ and $t\in \R$, from the following inequality
\begin{eqnarray*}
\|Q_{1}u(t)-Q_{\varepsilon}u(t)\|&\leq& C\int^{t}_{t-\varepsilon}\|T(t-s)F(s,u(s),u_{s})\|ds\\[8pt] &\leq& C\int^{t}_{t-\varepsilon}\|T(t-s)F(s,u(s),u_{s})\|ds\\[8pt]
&\leq& CMM'\varepsilon, \end{eqnarray*}
which implies that the set $(Q_{1}\overline{\Omega}_{R})(t)$ is relatively compact in $X$ for all $t\in \R$.

Thus, the Arzela-Ascoli theorem guarantees that $Q_{1}$ is a compact operator.

Secondly, we prove that  $Q_{2}$ is a contraction. Let $u,v\in \overline{\Omega}_{R}$, by the condition (H2), we have
\begin{eqnarray*}
&&\|Q_{2}u(t)-Q_{2}v(t)\|\\[8pt]
&\leq&\|(I-T(\omega))^{-1}\|\cdot\Big\|\sum_{t-\omega\leq t_{k}<t}T(t-t_{k})(I_{k}(u(t_{k}))-I_{k}(v(t_{k}))\Big\|\\[8pt]
&\leq& \frac{1}{1-e^{\nu_{0}\omega}}\cdot\sum_{k=1}^{p}M e^{\nu_{0}(t-t_{k})}a_{k}\|u(t_{k})-v(t_{k})\|\\[8pt]
&\leq& -\frac{M}{\nu_{0}\omega}\sum_{k=1}^{p}a_{k}\|u-v\|_{PC},
\end{eqnarray*}
therefore,
$$\|Q_{2}u-Q_{2}v\|_{PC}\leq -\frac{M}{\nu_{0}\omega}\sum_{k=1}^{p}a_{k}\|u-v\|_{PC}.\eqno(3.10)$$
From the condition (H3), we can deduce $Q_{2}$ is a contraction.

Therefore, by the famous Sadovskii fixed point theorem \cite{Sadovskii1967}, we know that $Q$ has a fixed point $u\in\overline{\Omega}_{R}$, that is, Eq. (1.1) has a $\omega$-periodic mild solution. The proof is completed. \hfill$\Box$

\vskip3mm
Furthermore, we assume that $F$ satisfies Lipschitz condition, namely,
\emph{\vskip1mm
\noindent (H1$'$) there are positive constants $c_{1},c_{2}$, such that for every $t\in\R, x_{0},x_{1}\in X$ and $\phi,\psi\in PC([-r,0],X)$
$$\|F(t,x,\phi)-F(t,y,\psi)\|\leq c_{1}\|x-y\|+c_{2}\|\phi-\psi\|_{Pr}, $$}
then we can obtain the following result.
\vskip3mm
\noindent\textbf{Theorem 3.2.}\emph{ Let $X$ be a Banach space, $-A$ generates an exponentially stable compact semigroup $T(t)(t\geq0)$ in $X$. Assume that  $F:\R\times X\times PC([-r, 0], X)\rightarrow X$ is continuous and $F(t,\cdot,\cdot)$ is $\omega$-periodic in $t$, $I_{k}\in C(X,X)(k\in \mathbb{Z})$.  If the conditions (H1 $'$),(H2) and (H3) hold, then Eq. (1.1) has unique $\omega$-periodic mild solution $u$.}
\vskip3mm
\noindent\textbf{Proof } From (H1$'$) we easily see that (H1) holds. In fact, for any $t\in\R$, $x\in X$ and $\phi\in PC([-r,0],X)$, by the condition (H1$'$),
\begin{eqnarray*}\|F(t,x,\phi)\|&\leq&\|F(t,x,\phi)-F(t,\theta,\theta)\|+\|F(t,\theta,\theta)\|\\[8pt]
&\leq&c_{1}\|x\|+c_{2}\|\phi\|_{Pr}+\|F(t,\theta,\theta)\|.\end{eqnarray*}
From the continuity and periodicity of $F$, we can choose $c_{0}=\max_{t\in[0,\omega]}\|F(t,\theta,\theta)\|$, thus, the condition (H1) holds. Hence by Theorem 3.1, Eq.(1.1) has $\omega$-periodic mild solutions.

Let $u,v\in PC_{\omega}(\R,X)$ be the $\omega$-periodic mild solutions of Eq.(1.1), then they are the fixed points of the operator $Q$ which is defined by (3.1). Hence,
\begin{eqnarray*}
&&\|u(t)-v(t)\|=\|Qu(t)-Qv(t)\|\\[8pt]
&\leq&\Big\|(I-T(\omega))^{-1}\int_{t-\omega}^{t}T(t-s)\Big(F(s,u(s),u_{s})-F(s,v(s),v_{s})\Big)ds\Big\|\\[8pt]
&&+\Big\|(I-T(\omega))^{-1}\sum_{t-\omega\leq t_{k}<t}T(t-t_{k})\Big(I_{k}(u(t_{k}))-I_{k}(v(t_{k}))\Big)\Big\|\\
&\leq&  \frac{1}{1-e^{\nu_{0}\omega}}\cdot\int_{t-\omega}^{t}Me^{\nu_{0}(t-s)}\cdot \|F(s,u(s),u_{s})-F(s,v(s),v_{s}))\|ds\\[8pt]
&&+\frac{1}{1-e^{\nu_{0}\omega}}\sum_{t-\omega\leq t_{k}<t}Me^{\nu_{0}(t-t_{k})}\cdot\|I_{k}(u(t_{k}))-I_{k}(v(t_{k}))\|\\[8pt]
&\leq& -\frac{M}{\nu_{0}}\int_{t-\omega}^{t}c_{1}\|u(s)-v(s)\|+c_{2}\|u_{s}-v_{s}\|_{Pr}ds-\frac{M}{\nu_{0}\omega}\sum_{k=1}^{p}a_{k}\|u(t_{k})-v(t_{k})\|\\[8pt]
&\leq&\Big(-\frac{M}{\nu_{0}}(c_{1}+c_{2})-\frac{M}{\nu_{0}\omega}\sum_{k=1}^{p}a_{k}\Big)\|u-p\|_{PC},
\end{eqnarray*}
which implies that $$\|u-v\|_{PC}=\|Qu-Qv\|_{PC}\leq \Big(-\frac{M}{\nu_{0}}(c_{1}+c_{2})-\frac{M}{\nu_{0}\omega}\sum_{k=1}^{p}a_{k}\Big)\|u-v\|_{PC}.$$ From this and the condition (H3), it follows that $u_{2}=u_{1}$. Thus, Eq.(1.1) has only one $\omega$-periodic mild solution.   \hfill$\Box$

\section{ The Asymptotic Stability}

\vskip3mm
In order prove the asymptotic stability of $\omega$-periodic solutions for  Eq. (1.1), we need discuss
 the existence and uniqueness of the following initial value problem
 $$\left\{\begin{array}{ll}
u'(t)+Au(t)=F(t,u(t),u_{t}),\ t\geq0,\ t\neq t_{i},\\[8pt]
 \Delta u(t_{i})=I_{i}(u(t_{i})),\ \qquad i=1,2,\cdots,\\[8pt]
 u_{0}=\varphi,
 \end{array} \right.\eqno(4.1)$$
where $F:J\times X\times PC([-r, 0], X)\rightarrow X$ is continuous and $\varphi\in PC([-r, 0], X)$.

Define
\begin{eqnarray*}PC([-r,\infty),X)&:=&\Big\{u:[-r,\infty)\to X | u\ \mathrm{is\ continuous\ at}\ t\in [-r,\infty)\setminus \{t_{j}\},\\[8pt]
&&u\ \mathrm{is\ continuous\ from\ left\ and\ has\ right\ limits\ at}\ t_{j}, j\in \N\Big\}.\end{eqnarray*} If there exists $u\in PC([-r,\infty),X)$ satisfying
   $u(t)=\varphi (t)$ for $-r\leq t\leq0$ and
  $$u(t)=T(t)u(0)+\int^{t}_{0}T(t-s)F(s,u(s),u_{s})ds+\sum_{0<t_{k}<t}T(t-t_{k})I_{k}(u(t_{k})),\ \ \  t\geq0, \eqno(4.2)$$
 then $u$ is
called a mild solution of the initial value problem (4.1).
 \vskip3mm

In order to obtain the results about asymptotic stability, we need  the following
integral inequality of Gronwall-Bellman type with delay and impulsive.
\vskip3mm
\noindent\textbf{Lemma 4.1.}\emph{ Assume that $\phi\in PC([-r,\infty),J)$, and
there exist constants
$\alpha_1\geq0, \alpha_2\geq0$ and $\beta_{i}\geq0(i=1,2\cdots,)$ such that for every $t\geq0$, $\phi$ satisfy the integral inequality
$$\phi(t)\leq\phi(0)+\alpha_{1}\int_{0}^{t}\phi(s)ds+\alpha_{2}\int_{0}^{t}\sup_{\tau\in[-r,0]}\phi(s+\tau)ds+\sum_{0<t_{k}<t}\beta_{k}\phi(t_{k}).\eqno(4.3)$$
 Then $\phi(t)\leq\|\phi\|_{Pr}\cdot\prod\limits_{0<t_{k}<t}(1+\beta_{k})e^{(\alpha_{1}+\alpha_{2})t}$ for every $t\geq0$.}
\vskip3mm
\noindent\textbf{Proof} Define a function $\psi:[-r,\infty)\to \R$ as following
 $$\psi(t)=\sup_{s\in[-r,t]}\phi(s),\ \ \ t\in[-r,\infty).$$
 Then $\psi\in PC([-r,\infty),J)$ and $\phi(t)\leq\psi(t)$ for $t\in[-r,\infty)$.
 Similar to the proof  of \cite[Lemma 4.1]{Li2011}, we can get the following inequality
 $$\psi(t)\leq \|\phi\|_{Pr}+(\alpha_{1}+\alpha_{2})\int_{0}^{t}\psi(s)ds+\sum_{0<t_{k}<t}\beta_{k}\psi(t_{k}),\ \qquad t\geq0,$$
 holds. By \cite[Lemma 1, p.12]{Samoilenko}, $\psi(t)\leq \|\phi\|_{Pr}\prod_{0<t_{k}<t}(1+\beta_{k})e^{(\alpha_{1}+\alpha_{2})t}$ for every $t\geq0$. Therefore,
 $\phi(t)\leq \|\phi\|_{Pr}\prod_{0<t_{k}<t}(1+\beta_{k})e^{(\alpha_{1}+\alpha_{2})t}$ for every $t\geq0$.  \hfill$\Box$

\vskip3mm
For the initial value problem (4.1), we have the following  result.
\vskip3mm
\noindent\textbf{Theorem 4.1.} \emph{ Let $X$ be a Banach space, $-A$ generates a $C_{0}$-semigroup $T(t)(t\geq0)$ in $X$. Assume that $F:J\times X\times PC([-r, 0], X)\rightarrow X$ is continuous, $I_{k}\in C(X,X)(k=1,2,\cdots)$, and $\varphi\in PC([-r,0],X)$. If the conditions (H1 $'$) and (H2) hold, then the initial value problem (4.1) has a
 unique  mild solution $u\in PC([-r,\infty),X)$.}
\vskip3mm
\noindent\textbf{Proof} For  $t\in[-r,t_{1}]$, the initial value problem (4.1)  is in the following form:
 $$\left\{\begin{array}{ll}
u'(t)+Au(t)=F(t,u(t),u_{t}),\ t\in[0,t_{1}],\\[8pt]
 u(t)=\varphi(t),\ \ \ t\in[-r,0].
 \end{array} \right.\eqno(4.3)$$
 Write
\begin{eqnarray*}PC([-r,t_{1}],X)&:=&\Big\{u |u:[-r,t_{1}]\to X\ \mathrm{ with}\ u|_{[0,t_{1}]}\in C([0,t_{1}],X)\\
 &&\mathrm{and }\ u|_{[-r,0]}\in PC([-r,0],X)\Big\},\end{eqnarray*}
 then $PC([-r,t_{1}],X)$ is a Banach space under the norm
 $$\|u\|_{1}=\sup_{t\in [0,t_{1}]}\|u(t)\|+\|u|_{[-r,0]}\|_{Pr}.$$
 For any $\varphi\in PC([-r,0],X)$, let $PC_{\varphi}([-r,t_{1}],X)=\{u\in PC([-r,t_{1}],X)|\ u|_{[-r,0]}=\varphi \}$, then $PC_{\varphi}([-r,t_{1}],X)$ is  a closed convex subset of $PC([-r,t_{1}],X)$. For each $\varphi\in PC([-r,0],X)$ and $u\in PC_{\varphi}([-r,t_{1}],X)$,

Define an operator as following
$$(\widetilde{Q}u)(t)=\left\{\begin{array}{ll}T(t)u(0)+\int^{t}_{0}T(t-s)F(s,u(s),u_{s})ds,\ \ \  t\in [0,t_{1}],\\[8pt]
\varphi(t), \ \qquad \qquad\qquad\qquad t\in[-r,0].\end{array} \right. \eqno(4.4)$$
 It is easy to see that $\widetilde{Q}$ is well defined, $\widetilde{Q}u\in PC_{\varphi}([-r,t_{1}],X)$, and the mild solution of Eq.(4.3) for $\phi$ is equivalent to the fixed point
of $\widetilde{Q}$  in $PC_{\varphi}([-r,t_{1}],X)$.

Now, we prove that $\widetilde{Q}$ has a fixed point in $PC_{\varphi}([-r,t_{1}],X)$. From the condition (H1 $'$) and (4.4), it follows that  for any $u, v\in PC_{\varphi}([-r,t_{1}],X)$, $n=1,2,\cdots$,
$$\|(\widetilde{Q}^{n}u)(t)-(\widetilde{Q}^{n}v)(t)\|\leq \frac{(M_{1}(c_{1}+c_{2})t)^{n}}{n!}\|u-v\|_{1},$$
where  $M_{1}$ is the bound of $\|T(t)\|$ on $[0,t_{1}]$.
By the contraction principle, one shows that
 $\widetilde{Q}$  has a unique fixed point $u_{1}$ in $PC_{\varphi}([-r,t_{1}],X)$,  which means the initial value problem (4.3) has a mild
solution and
 $$u_{1}(t)=\left\{\begin{array}{ll}T(t)u(0)+\int^{t}_{0}T(t-s)F(s,u(s),u_{s})ds,\\[8pt]
 \varphi(t),\ \ \ t\in [-r,0].
 \end{array} \right.\eqno(4.5)$$

For  $t\in[-r,t_{2}]$, the initial value problem (4.1)  is in the following form:
 $$\left\{\begin{array}{ll}
u'(t)+Au(t)=F(t,u(t),u_{t}),\ t\in(t_{1},t_{2}],\\[8pt]
 u(t^{+}_{1})=I_{1}(u(t_{1}))+u(t_{1}),\\[8pt]
 u(t)=u_{1}(t),\ \ \ \ t\in [-r,t_{1}].
 \end{array} \right.\eqno(4.6)$$
Similar to the proof of (4.3), we can prove that the initial value problem (4.6) has a mild solution $u_{2}\in PC_{\varphi}([-r,t_{2}],X)$
\begin{eqnarray*}u_{2}(t)&=&\left\{\begin{array}{ll}T(t-t_{1})u(t_{1}^{+})+\int^{t}_{t_{1}}T(t-s)F(s,u(s),u_{s})ds, \  t\in(t_{1},t_{2}],\\[8pt]
 u_{1}(t),\  t\in [-r,t_{1}],
\end{array} \right.\\[8pt]
 &=&\left\{\begin{array}{ll}T(t)u(0)+\int^{t}_{0}T(t-s)F(s,u(s),u_{s})ds+T(t-t_{1})I_{1}(u(t_{1})),\ t\in [0,t_{2}], \\[8pt]
 \varphi(t),\ \ \ t\in [-r,0].\end{array} \right.\end{eqnarray*}
Doing this interval by interval, we obtain that there exists $u\in PC_{\varphi}([-r,\infty),X)$ satisfying
   $u(t)=\varphi (t)$ for $-r\leq t\leq0$ and
  $$u(t)=T(t){u}(0)+\int^{t}_{0}T(t-s)F(s,{u}(s),u_{s})ds+\sum_{0<t_{k}<t}T(t-t_{k})I_{k}(u(t_{k})),\ \ \  t\geq0,\eqno(4.7)$$
which is a mild solution of the
the initial value problem (4.1).

Next, we show the uniqueness. Let  $u,v\in PC([-r,\infty),X)$
be the mild  solutions of the initial value problem (4.1),
hence they satisfy  the initial value condition $u(t)=v(t)=\varphi(t)$ for $-r\leq t\leq0$ and (4.2).
 By the condition (H1 $'$) and (H2),
  for every $t\geq0$, one has
\begin{eqnarray*}
&&\|u(t)-v(t)\|\\[8pt]
&\leq&\Big\|\int_{0}^{t}T(t-s)\Big(F(s,u(s),u_{s})-F(s,v(s),v_{s})\Big)ds\Big\|\\
&&\qquad\qquad\qquad\qquad\qquad+\Big\|\sum_{0<t_{k}<t}T(t-t_{k})(I_{k}(u(t_{k}))-I_{k}(v(t_{k})))\Big\|\\
&\leq&\int_{0}^{t}\|T(t-s)\|\cdot\|F(s,u(s),u_{s})-F(s,v(s),v_{s})\|ds\\
&&\qquad\qquad\qquad\qquad\qquad+\sum_{0<t_{k}<t}\|T(t-t_{k})\|\cdot\|I_{k}(u(t_{k}))-I_{k}(v(t_{k}))\|\\
&\leq&M\int_{0}^{t}c_{1}\|u(s)-v(s)\|+c_{2}\|u_{s}-v_{s}\|_{Pr}ds\\
&&\qquad\qquad\qquad\qquad\qquad+M\sum_{0<t_{k}<t}a_{k}\|u(t_{k})-v(t_{k})\|\\
&\leq&M\int_{0}^{t}c_{1}\|u(s)-v(s)\|+c_{2}\sup\limits_{\tau\in[-r,0]}\|u(s+\tau)-v(s+\tau)\|ds\\
&&\qquad\qquad\qquad\qquad\qquad+M\sum_{0<t_{k}<t}a_{k}\|u(t_{k})-v(t_{k})\|.
\end{eqnarray*}

From Lemma 4.1, it follows that $\|u(t)-v(t)\|=0$ for every $t\geq0$.
Hence, $u\equiv v$. This completes the proof of Theorem 4.1.
\hfill$\Box$
\vskip3mm

\vskip3mm
\noindent\textbf{Theorem 4.2.}\emph{ Let $X$ be a Banach space, $-A$ generates an exponentially stable compact semigroup $T(t)(t\geq0)$ in $X$. Assume that  $F:\R\times X\times PC([-r, 0], X)\rightarrow X$ is continuous and $F(t,\cdot,\cdot)$ is $\omega$-periodic in $t$, $I_{k}\in C(X,X)(k\in \mathbb{Z})$.  If the conditions (H1 $'$),(H2) and
\vskip2mm
\noindent (H3 $'$)  $(c_{1}+c_{2}e^{-\nu_{0}r})+\frac{1}{\omega}\sum\limits_{k=1}^{p}a_{k}<\frac{|\nu_{0}|}{M}$,
\vskip2mm
\noindent hold, then  the unique  $\omega$-periodic
mild solution of the periodic problem (1.1) is globally asymptotically stable.}
\vskip3mm
\noindent\textbf{Proof } From the condition (H3 $'$), it follows that the condition (H3) holds. By Theorem 3.2, the periodic problem (1.1) has a unique $\omega$-periodic mild solution $u^{*}\in PC_{\omega}(\R,X)$. For any $\phi\in PC([-r,0],X)$,
the initial value problem (4.1) has a unique global mild solution $u=u(t,\phi)\in PC([-r,\infty),X)$ by Theorem 4.1.

By the semigroup representation of the solutions, $u^{*}$ and $u$
satisfy the integral equation (4.2). Thus, by (4.2) and
condition (H1 $'$), (H2), for any $t\geq0$, we have
 \begin{eqnarray*}&&\|u(t)-u^{*}(t)\|\\[8pt]
 &\leq&\Big\|T(t)(u(0)-u^{*}(0))\Big\|
 +\Big\|\int^{t}_{0}T(t-s)(c_{1}(u(s)-u^{*}(s))
 +c_{2}(u_{s}-u^{*}_{s}))ds\Big\|\\[8pt]
 &&+\Big\|\sum_{0<t_{k}<t}T(t-t_{k})\Big(I_{k}(u(t_{k}))-I_{k}(u^{*}(t_{k}))\Big)\Big\|\\[8pt]
 &\leq& Me^{\nu_{0} t}\|u(0)-u^{*}(0)\|+\int^{t}_{0}Me^{\nu_{0}(t-s)}c_{1}\|u(s)-u^{*}(s)\|ds\\[8pt]
 &&+\int^{t}_{0}Me^{\nu_{0}(t-s)}c_{2}\sup_{\tau\in[-r,0]}\|u(s+\tau)-u^{*}(s+\tau)\|ds\\[8pt]
 &&+\sum_{0<t_{k}<t}Me^{\nu_{0}(t-t_{k})}a_{k}\|u(t_{k})-u^{*}(t_{k})\|ds\\[8pt]
 &\leq& Me^{\nu_{0} t}\|u(0)-u^{*}(0)\|+e^{\nu_{0}t}\int^{t}_{0}Mc_{1}e^{-\nu_{0}s}\|u(s)-u^{*}(s)\|\\[8pt]
 &&+e^{\nu_{0}t}\int^{t}_{0}Mc_{2}e^{-\nu_{0}r}\sup\limits_{\tau\in[-r,0]}e^{-\nu_{0}(s+\tau)}\|u(s+\tau)-u^{*}(s+\tau)\|ds\\[8pt]
 &&+e^{\nu_{0}t}\sum_{0<t_{k}<t}M a_{k}e^{-\nu_{0}t_{k}}\|u(t_{k})-u^{*}(t_{k})\|.\end{eqnarray*}
Then
\begin{eqnarray*}e^{-\nu_{0}t}\|u(t)-u^{*}(t)\|
 &\leq& M\|u(0)-u^{*}(0)\|+\int^{t}_{0}M c_{1}e^{-\nu_{0}s}\|u(s)-u^{*}(s)\|\\[8pt]
 &&+\int^{t}_{0}M c_{2}e^{-\nu_{0}r}\sup\limits_{\tau\in[-r,0]}e^{-\nu_{0}(s+\tau)}\|u(s+\tau)-u^{*}(s+\tau)\|ds\\[8pt]
 &&+\sum_{0<t_{k}<t}M a_{k}e^{-\nu_{0}t_{k}}\|u(t_{k})-u^{*}(t_{k})\|,\end{eqnarray*}
for $t\in[-r,\infty)$, let $\phi(t)=e^{-\nu_{0}t}\|u(t)-u^{*}(t)\|$,
 one can obtain
 $$\phi(t)\leq M\phi(0)+M c_{1}\int^{t}_{0}\phi(s)+M c_{2}e^{-\nu_{0}r}\int^{t}_{0}\sup\limits_{\tau\in[-r,0]}\phi(s+\tau)ds+\sum_{0<t_{k}<t}M a_{k}\phi(t_{k}),$$
 hence, from Lemma 4.1, it follows that
 $$e^{-\nu_{0}t}\|u(t)-u^{*}(t)\|=\phi(t)\leq C(\varphi)\prod_{0<t_{k}<t}(1+M a_{k})e^{M(c_{1}+c_{2}e^{-\nu_{0}r})t},\eqno(4.8)$$
where \ $C(\varphi)=\sup\limits_{s\in[-\tau,0]}\{e^{-\nu_{0}s}\|\varphi(s)-u^{*}(s)\|_{Pr}\}$.
Set $k=np+q(q=1,2\cdots,p-1, n=0,1,2,\cdots)$, by the periodicity of $a_{k}$, one can obtain
$$\sum_{0<t_{k}<t}\ln(1+M a_{k})\leq (n+1)\sum_{k=1}^{p}\ln(1+M a_{k}),$$
thus,
\begin{eqnarray*}\lim_{t\to\infty}\frac{\sum\limits_{0<t_{k}<t}\ln(1+M a_{k})}{t}&\leq&\lim_{n\to\infty} \frac{(n+1)\sum\limits_{k=1}^{p}\ln(1+M a_{k})}{n\omega}\\[8pt]
&=&\frac{1}{\omega}\sum\limits_{k=1}^{p}\ln(1+M a_{k})\\[8pt]&\leq&\frac{M }{\omega}\sum\limits_{k=1}^{p}a_{k}
.\end{eqnarray*}
Therefore, from the assumption (H3), it follows that
$$\sigma:=-\nu_{0}-\lim_{t\to\infty}\frac{\sum_{0<t_{k}<t}\ln(1+M a_{k})}{t}-M(c_{1}+c_{2}e^{-\nu_{0}r})>0,\eqno(4.9)$$
Combining (4.8) with (4.9), we can obtain
$$\|u(t)-u^{*}(t)\|\leq C(\varphi)e^{-\sigma t}\to 0\ \ \ (t\to \infty).\eqno(4.10)$$
Thus, the $\omega$-periodic solution $u^{*}$ is
globally asymptotically stable and it exponentially
attracts every mild solution of the initial value problem.
 This completes the proof of Theorem 4.2.
\hfill$\Box$

\section{Application}
In this section, we present one example, which does not aim at generality, but indicates
how our abstract results can be applied to  concrete problems.

Let $\overline{\Omega}\in\R^{n}$ be a bounded domain with
a $C^{2}$-boundary $\partial\Omega$ for $n\in N$.
Let $\nabla ^{2}$ is a Laplace operator, and $\lambda_{1}$ is the smallest eigenvalue of operator $-\nabla ^{2}$ under the Dirichlet boundary condition $u|_{\partial\Omega}=0$. It is well known (\cite[Theorem 1.16]{Amann76},) that $\lambda_{1}>0$.

Under the above assumptions,  we discuss the existence,
uniqueness and asymptotic stability of time
$2\pi$-periodic solutions of the semilinear parabolic
boundary value problem
$$\left\{\begin{array}{ll}
\frac{\partial u}{\partial t}-\nabla ^{2}u=\frac{\lambda_{1}}{4}\sin t\cdot u(x,t)+\int_{t-r}^{t}e^{\frac{4}{\lambda_{1}}(s-t) }\cdot u(x,s)ds,\ x\in \Omega,  t\neq t_{k},\\[10pt]
\triangle u(x,t_{k})=\frac{\lambda_{1}\pi}{p}(e^{\sin u(x,t_{k})}-1),\ \ \ \ x\in \Omega,\ t_{k}=\frac{2k-1}{p}\pi, k\in\Z,\\[10pt]
u|_{\partial\Omega}=0,
 \end{array} \right.\eqno (5.1)$$
where $p>0$ is a integer, $r>0$ is a real number.

Let\ $X=L^{2}(\Omega)$ with the norm $\|\cdot\|_{2}$,
then\ $X$ is a Banach space.
Define an operator\ $A:D(A)\subset X\rightarrow X$ by:
$$D(A)=  W^{2,2}(\Omega)\cap W_{0}^{1,2}(\Omega), \quad Au=-\nabla ^{2}u.\eqno(5.2)$$
From \cite{Amann78}, we know that $-A$ is a selfadjoint operator in $X$,
and generates an exponentially stable analytic semigroup $T_{p}(t)(t\geq0)$, which is contractive in $X$. From the operator $A$ has compact resolvent in $L^{2}(\Omega)$,
$T_{p}(t)(t\geq0 )$ is a compact semigroup (see \cite{Pazy}), which implies that the growth exponent of the semigroup $T(t)(t\geq0)$ satisfies $\nu_{0}=-\lambda_{1}$. Therefore, for every $t>0$, $\|T(t)\|_{2}\leq M:=1$ and $\|(I-T(2\pi))\|\leq \frac{1}{1-e^{-2\lambda_{1}\pi}}$.

Now, we define $u(t)=u(\cdot,t)$,  $I_{k}(u(t_{k}))=(e^{\sin u(\cdot,t_{k})}-1)$, and since for $u\in PC_{2\pi}(\R,X)$,
$$\int_{t-r}^{t}e^{\frac{4}{\lambda_{1}}(s-t) }\cdot u(s)ds=\int_{-r}^{0}e^{\frac{4}{\lambda_{1}}s }\cdot u(t+s)ds=\int_{-r}^{0}e^{\frac{4}{\lambda_{1}}s }\cdot u_{t}(s)ds,$$
we define $F:\R\times X\times PC([-r,0],X)\to X$ by
$$F(t,\xi,\phi)=\frac{\lambda_{1}}{4}\sin t\cdot \xi+\int_{-r}^{0}e^{\frac{4}{\lambda_{1}} s}\phi(s) ds, \eqno(5.3)$$
thus, it is east to see that $t_{k+p}=t_{k}+2\pi$,  $I_k : X\rightarrow X(k\in\Z)$ are continuous functions satisfying $I_{k+p}=I_{k}$ and $F:\R\times X\times PC([-r, 0], X)\rightarrow X$  is continuous function which is $2\pi$-periodic in $t$. Hence, the impulsive and delay parabolic boundary value problem
(5.1) can be reformulated as the abstract evolution equation (1.1) in $X$.

 From the definition of $F$ and $I_{k}$, for every $t\in\R,\xi_{1},\xi_{2}\in X$ and $\phi_{1},\phi_{2}\in PC([-r,0],X)$, we have
\begin{eqnarray*}\|F(t,\xi_{1},\phi_{1})-F(t,\xi_{2},\phi_{2})\|_{2}&\leq& \frac{\lambda_{1}}{4}\|\xi_{1}-\xi_{2}\|_{2}+
\frac{\lambda_{1}}{4}(1-e^{-\frac{4r}{\lambda_{1}}})\|\phi_{1}-\phi_{2}\|_{Pr},\\[8pt]
\|I_{k}(\xi_{1})-I_{k}(\xi_{2})\|_{2}&\leq& \frac{\lambda_{1}\pi}{p}\|\xi_{1}-\xi_{2}\|_{2},\ \ \ k\in\Z,
\end{eqnarray*}
which implies that the conditions (H1 $'$), (H2) and (H3) hold. Thus, by the Theorem 3,2, the parabolic boundary value problem (5.1) has only one time $2\pi$-periodic mild solution.
Moreover, if $0<r<\frac{\lambda_{1}\ln 4}{\lambda_{1}^{2}+4}$, then we can deduce that the condition (H3 $'$) holds.  From the Theorem 4.2, one can see that the unique $2\pi$-periodic mild solution of problem (5.1) is globally asymptotically stable.

\bibliographystyle{abbrv}

\end{document}